\documentclass[12pt]{article}
\usepackage{amssymb,amsmath}

\usepackage{latexsym,amsfonts,amssymb,amsthm,amscd,amsmath,graphicx}
\usepackage{amsmath}
\usepackage{amsfonts}
\usepackage{amssymb}
\usepackage{graphicx}%
\usepackage{booktabs}

\begin{document}
\title{The Adaptive Lasso with Weighted IPW Factor and KRR Penalty Function}

\author{
  Zeen Song\thanks{E-mail: zeensong17@gmail.com} \\[0.5em]
  \small School of Mathematics, Statistics and Mechanics,Beijing University of Technology
}

\maketitle
\date{}

\newtheorem{theorem}{Theorem}[section]
\newtheorem{definition}{Definition}[section]
\newtheorem{lemma}{Lemma}[section]
\newtheorem{proposition}{Proposition}[section]
\newtheorem{corollary}{Corollary}[section]
\newtheorem{remark}{Remark}[section]
\renewcommand{\theequation}{\thesection.\arabic{equation}}
\catcode`@=11 \@addtoreset{equation}{section} \catcode`@=12

\begin{abstract}
In this paper, we first consider a Lasso estimator incorporating a weighted IPW factor, and obtain certain conditions for the Lasso variable selection to be consistent or inconsistent. Then we provide a modified version of the adaptive Lasso equipped with a weighted IPW factor and a KRR penalty function, establishing that the proposed adaptive Lasso possesses the Oracle property. Numerical simulations demonstrate that the proposed adaptive Lasso exhibits near-minimax optimality in performance.

\vspace{0.5em}
\end{abstract}

\section{Introduction}
\qquad Variable selection and high prediction accuracy are two fundamental goals in statistical learning. To attain the two aims, statisticians introduced a large number of predictors, meanwhile, they usually used stepwise deletion and subset selection. But the subset selection has several drawbacks such as lack of stability. Then a unified approach via penalized least squares was recommended, which contains ``Lasso-type" estimator and retains good features of both subset selection and ridge regression(see\cite{Tibshirani1996, KnightFu2000, FanLi2001, FanPeng2004, Meinshausen2004, Zou2006} and the references therein).

One of the consistent variance estimators is called imputation regression estimator which often appears in missing at random framework. Kernel ridge regression(KRR) is a general tool for nonparametric regression and inverse probability weighting (IPW) method is commonly used to analyze non-ignorable missing data(see\cite{Rubin1976, Cheng1994, WangKim2021, LiQinLiu2023} and the references therein).

Inspired by these references above, we discover a relevant predictive variable through defining a adaptive lasso with weighted IPW factor and KRR penalty function, which can be applied to high-dimensional nonparametric modeling.

Assume that $X_j=(x_{1j},...,x_{nj})^T, j=1,...,p$ are the linearly independent predictors and  $Y=(y_1,...,y_n)^T$ is the response vector, $p$ can be a fixed constant or a number keeps the relation to $n$ such as $p>>n$ or $=O(exp(n^{\theta})$ for some constant $\theta>0$. Let $\delta_i$ be the response indicator function of unit $i$ and $y_i$ is observed if and only if $\delta_i=1$.

Denote $$X=(X_1,...,X_p)=\left(
\begin{array}{ccc}
    \tilde{X}_{1p}\\
  ... \\
  \tilde{X}_{np}\\
 \end{array}\right)
 $$
the predictor matrix, $E(Y|X)=m(X)$ and its nonparametric estimator
\begin{align}
\hat{m}(X)=\frac{\sum_{i=1}^n\sum_{j=1}^p\delta_iK_h(x_{ij},X_j)y_i}{\sum_{i=1}^n\sum_{j=1}^p\delta_iK_h(x_{ij},X_j)},\label{425s1}
\end{align}
where $K_h(\cdot)$ is the kernel function with bandwidth $h$.

Suppose that

(C1). $m(x)=\beta^*_1X_1+...+\beta^*_pX_p$ and $|\mathcal{A}|=p_0<p$, where $\mathcal{A}=\{j:\beta^*_j\neq 0\}$.

(C2). $y_i=\tilde{X}_{ip}\beta^*+\epsilon_i, i=1,...,n$,  $\epsilon_1$,...,$\epsilon_n$ are independent identically distributed
 random variables with mean $0$ and variance $\sigma^2$.

 (C3). $\frac{1}{n}X^TAX\rightarrow B$ with $B$ is a positive definite matrix, where
  $$A=(a_{ij})_{n\times n}=diag \left((1-\frac{\delta_1}{\pi(\tilde{x}_1)})^2, (1-\frac{\delta_2}{\pi(\tilde{x}_2)})^2,..., (1-\frac{\delta_n}{\pi(\tilde{x}_n)})^2\right)_{n\times n}$$
  and
  $$
  B=\left(
       \begin{array}{cc}
         B_{11} & B_{12}\\
          B_{21} & B_{22} \\
       \end{array}
     \right),
  $$
where $B_{11}$ is a $p_0\times p_0$ matrix,
\begin{align}
 \pi(\tilde{x}_i)=\frac{1}{1+\exp(\alpha+\tilde{X}_{ip}\gamma_1+Y^T\gamma_2) }.\label{263281}
\end{align}
Here $\alpha$ is a constant, $\gamma_1$ is a known $p$-dimension vector and $\gamma_2$ is a known $n$-dimension vector.

Under the assumptions (C1)--(C3), we first consider the following lasso estimates
\begin{align}
\hat{\beta}^{(n)}&=arg\min_{\beta} \left[\sum_{i=1}^n (1-\frac{\delta_i}{\pi(\tilde{x}_i)})^2(y_i-\sum_{j=1}^p\beta_jx_{ij})^2+\lambda_n \sum_{j=1}^p|\beta_j|^{\gamma}\right]\nonumber\\
&=arg\min_{\beta} \left[(Y-\sum_{j=1}^pX_j\beta_j)^TA(Y-\sum_{j=1}^pX_j\beta_j)+\lambda_n \sum_{j=1}^p|\beta_j|^{\gamma}\right].\label{263231}
\end{align}

Unlike the classic Lasso-type estimator, we introduce a weighted IPW factor $\pi(\tilde{x}_i)$ in the first term, while we keep the second term ``$l_1$ --penalty" remain unchanged. First, we are concerned with the asymptotic property when $\lambda_n$ varies with $n$ as well as the conditions for the consistency of the lasso variable selection. Later, through a weighted KRR factor in the $l_1$ penalty, we define an adaptive Lasso  estimator and show that it enjoys the oracle properties.

The rest of this paper is organized as follows. In Section 2, we will deal with the conditions for the lasso
variable selection whether it maintains consistency or fails to be consistent. In Section 3, we will introduce a adaptive weighted lasso and prove that it enjoys the oracle properties. In Section 4, we will use a simulation study to compare our weighted lasso with the classic lasso and adaptive lasso. To be readable, we relegate the proofs of the Lemmas and Theorems to the Appendix.

\section{Conditions for Consistency of the Weighted IPW Lasso Variable Selection}

\qquad Without loss generality, denote $\mathcal{A}=\{1,2,...,p_0\}$ and $\mathcal{A}_n=\{j:\hat{\beta}^{(n)}_j\neq 0\}$.

First we state serval lemmas which consider the asymptotic behavior for lasso estimator when $\lambda_n$ varies with $n$.

\quad{\bf Lemma 1.} {\it If $\frac{\lambda_n}{n}\rightarrow \lambda_0\geq 0$, then
$\hat{\beta}^{(n)}\rightarrow_p arg\ min(Z)$, where
$$
Z(\beta)=(\beta-\beta^*)^TB(\beta-\beta^*)+\lambda_0\sum_{j=1}^p|\beta_j|^{\gamma}.
$$
 }

{\bf Lemma 2.} {\it Suppose that $\gamma\geq 1$. If $\frac{\lambda_n}{\sqrt{n}}\rightarrow \lambda_0\geq 0$ and $B$ is nonsingular,
then $\sqrt{n}(\hat{\beta}^{(n)}-\beta)\rightarrow_d arg\ min(V_1)$, where
$$
V_1(u)=-2u^Tw+u^TBu+\lambda_0\sum_{j=1}^pu_jsgn(\beta_j)|\beta_j|^{\gamma-1}\quad {\rm if }\ \gamma>1,
$$
$$
V_1(u)=-2u^Tw+u^TBu+\lambda_0\sum_{j=1}^p[u_jsgn(\beta_j)I(\beta_j\neq 0)+|u_j|I(\beta_j=0)]\quad {\rm if }\ \gamma=1,
$$
$w$ has a $N(0,\delta^2 B)$ distribution.
}

{\bf Lemma 3.} {\it Suppose that $\gamma<1$. If $\frac{\lambda_n}{n^{\frac{\gamma}{2}}}\rightarrow \lambda_0\geq 0$, then
$\sqrt{n}(\hat{\beta}^{(n)}-\beta)\rightarrow_d arg\ min(V_2)$, where
$$
V_2(u)=-2u^Tw+u^TBu+\lambda_0\sum_{j=1}^p|u_j|^{\gamma}I(\beta_j=0),
$$
$w$ has a $N(0,\delta^2 B)$ distribution.
}

{\bf Proposition 1.} {\it If $\frac{\lambda_n}{\sqrt{n}}\rightarrow \lambda_0\geq 0$, then
$\lim sup_n P(\mathcal{A}_n=\mathcal{A})\leq c<1$, where $c$ is a constant depending on the true model.
}

{\bf Lemma 4.} {\it If $\frac{\lambda_n}{n}\rightarrow 0$, $\frac{\lambda_n}{\sqrt{n}}\rightarrow \infty$, then
$\frac{n}{\lambda_n}(\hat{\beta}^{(n)}-\beta^*)\rightarrow_p arg\ min(V_3)$, where
$$
V_3(u)=u^TBu+\sum_{j=1}^p[u_jsgn(\beta^*_j)I(\beta^*_j\neq 0)+|u_j|I(\beta^*_j=0)].
$$
}

Now we give a necessary condition on the consistency of the lasso selection below.

{\bf Theorem 1(Necessary condition).} {\it
Assume that $\lim_n P(\mathcal{A}_n=\mathcal{A})=1$. Then there exists some sign vector $\mathbf{s}=(s_1,...,s_{p_0})^T$, $s_j=1$ or $-1$, such that
\begin{align}
|B_{21}B^{-1}_{11}\mathbf{s}|\leq 1\label{426s1}
\end{align}
in the sense of componentwise.
  }

The following corollary shows that  the condition (\ref{426s1}) is nontrivial because the lasso variable selection maybe inconsistent if it fails.

{\bf Corollary 1.} {\it Assume that $p_0=2m+1\ge 3$ and $p=p_0+1$. Choosing $\rho_1$ and $\rho_2$ satisfying $\frac{(n+1)\rho_2}{(n+1)-\rho_1}>1$, letting $B_{11}=(1-\rho_1)I+\rho_1J_1$, $B_{21}=(\rho_2,...,\rho_2)$ and $B_{22}=1$, where $J_1$ is the matrix with all of the elements equal to 1, it is easy to verify the condition (\ref{426s1}) cannot be satisfied. Therefore, the lasso variable selection is inconsistent. }

The results above show that the weighted IPW lasso cannot be an oracle procedure. Naturally, we are interested in which lasso can enjoy the oracle property.

\section{The Adaptive Lasso with Weighted KRR Penalty Function}

\qquad In this section,  we define the following adaptive lasso with weighted KRR penalty function
\begin{align}
\hat{\beta}^{*(n)}=arg\min_{\beta} \left[\sum_{i=1}^n (1-\frac{\delta_i}{\pi(\tilde{x}_i)})^2(y_i-\sum_{j=1}^p\beta_jx_{ij})^2+\lambda_n \sum_{j=1}^p\hat{w}_j|\beta_j|\right],\label{26328bc2}
\end{align}
where $\hat{w}=|\hat{m}(X)|\frac{1}{|\hat{\beta}|^{\gamma}}(1,...,1)^T$. (\ref{26328bc2}) can be looked as a remedy of (\ref{263231}) with a simple modification and it is an effective one. In fact, we will show that the adaptive lasso (\ref{26328bc2}) enjoys the following oracle property:

{\bf Theorem 2(Oracle properties)} {\it Assume that $\frac{\lambda_n}{\sqrt{n}}\rightarrow 0$ and $\lambda_nn^{\frac{(\gamma-1)}{2}}\rightarrow \infty$.
Then

1.  $$
\lim_n P(\mathcal{A}^*_n=\mathcal{A})=1,
$$

2. $$
\sqrt{n}(\hat{\beta}^{*(n)}_{\mathcal{A}}-\beta^*_{\mathcal{A}})\rightarrow_d N(0,\sigma^2B^{-1}_{11}).
$$

Here $\mathcal{A}^*_n=\{j:\hat{\beta}^{*(n)}_j\neq 0\}$. That is, the adaptive lasso estimates satisfy consistency in variable selection and asymptotic normality.
}

\section{Numerical Experiments}

\qquad In this section, we present simulation examples and compare the MSEs and RPEs of the Lasso and adaptive Lasso estimators across 200, 400, and 800 Monte Carlo replications.

\begin{center}
\includegraphics[height=6cm]{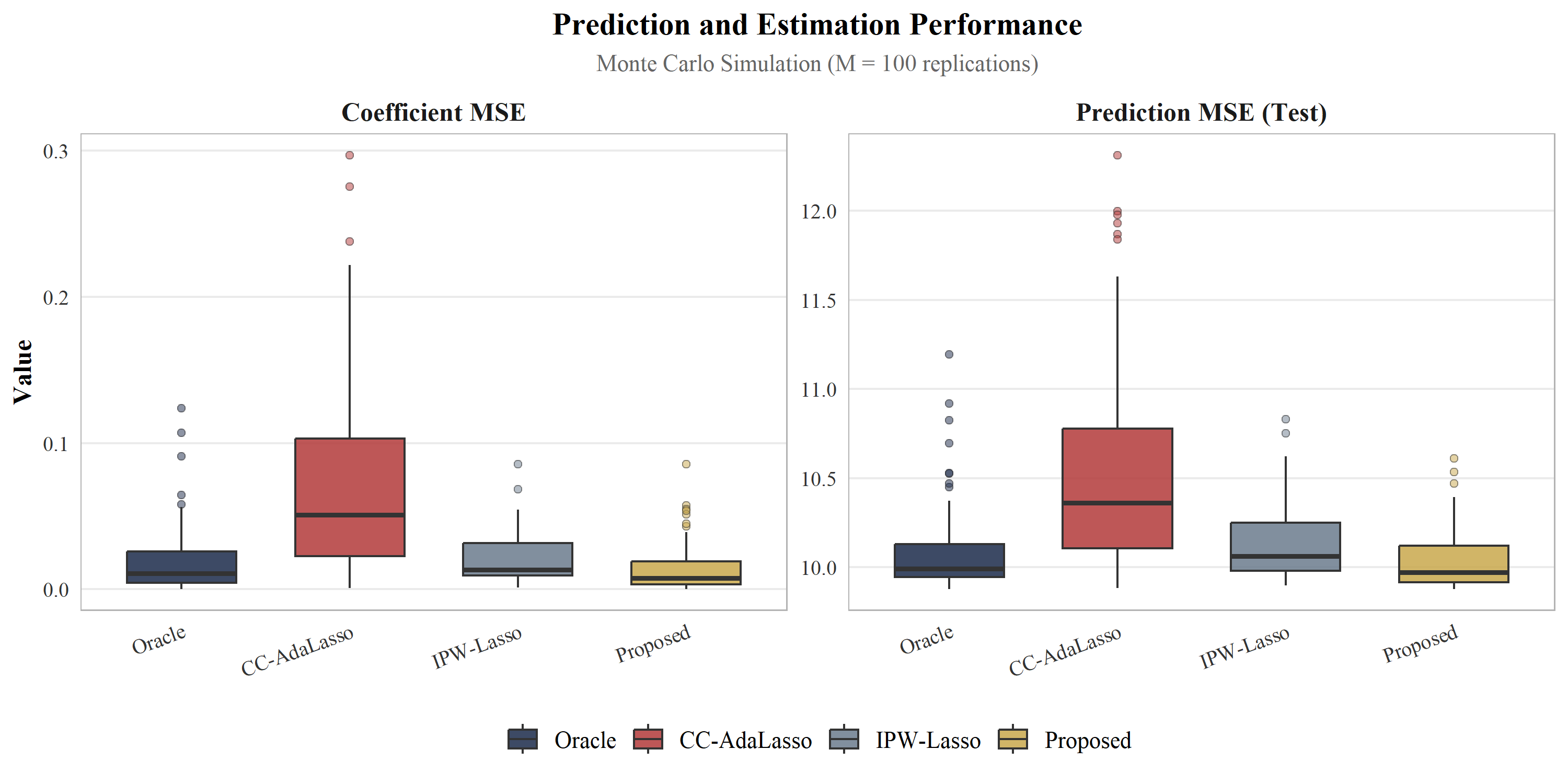}
\end{center}
\begin{center}
\underline{Figure 1.n=200 Regression Model }
\end{center}

\begin{center}
\includegraphics[height=6cm]{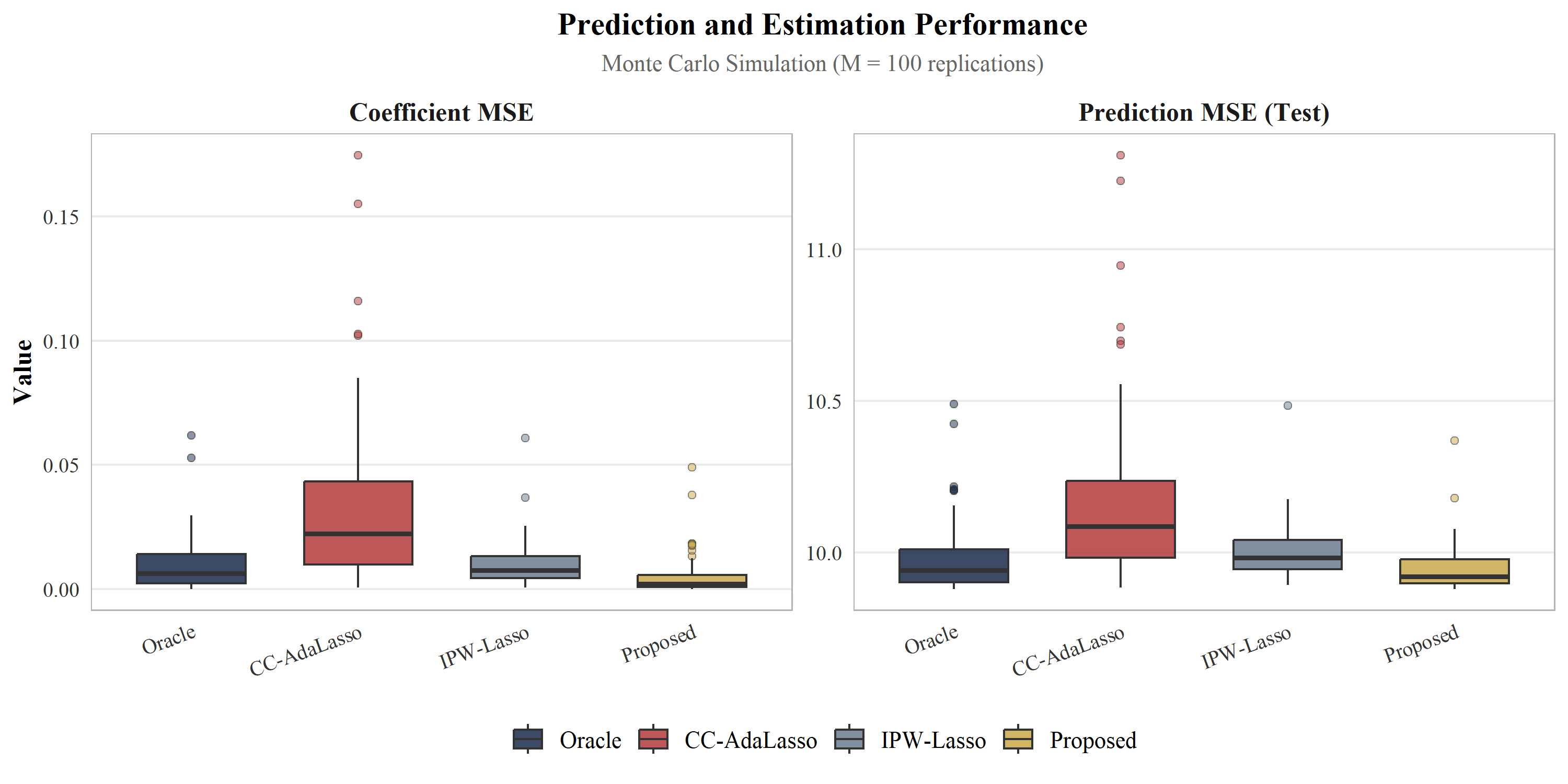}
\end{center}
\begin{center}
\underline{Figure 2.n=400 Regression Model }
\end{center}

\begin{center}
\includegraphics[height=6cm]{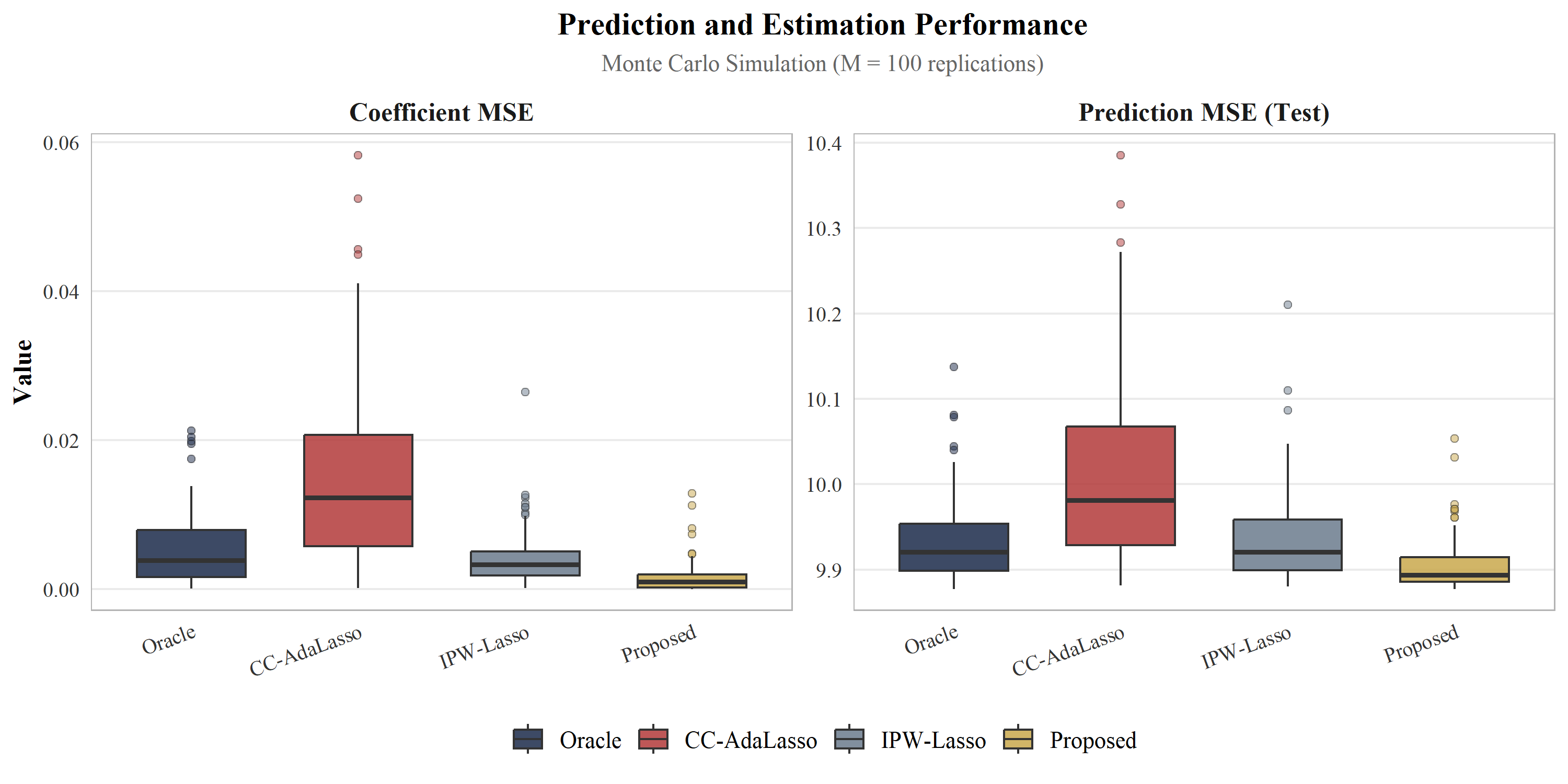}
\end{center}
\begin{center}
\underline{Figure 3.n=800 Regression Model }
\end{center}

The simulation results demonstrate that the proposed IPW adaptive Lasso consistently achieves superior performance across all candidate sample sizes.

We also report the relative prediction errors (RPEs) under different settings:

\begin{center}
\includegraphics[height=6cm]{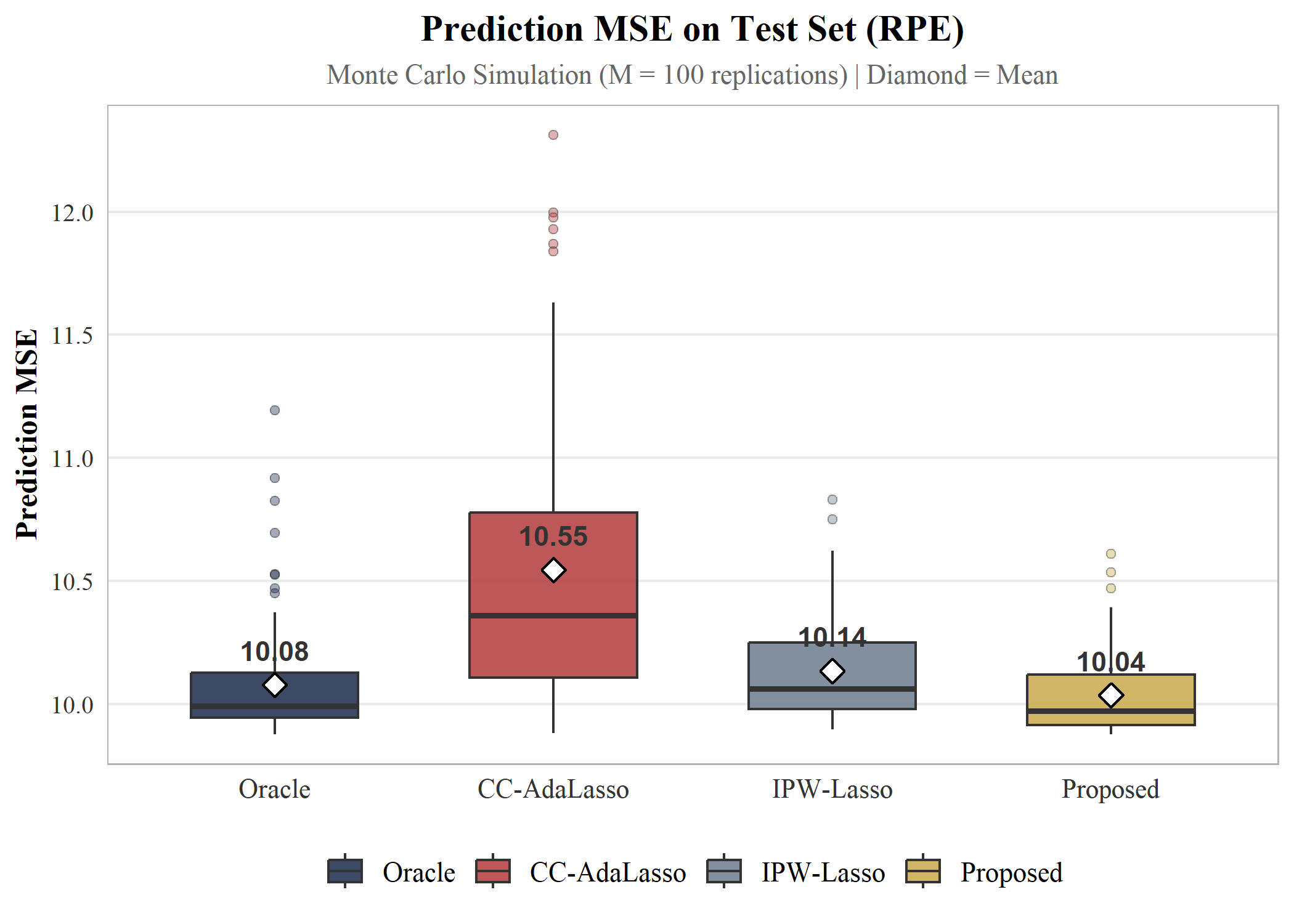}
\end{center}
\begin{center}
\underline{Figure 4.n=200  }
\end{center}

\begin{center}
\includegraphics[height=6cm]{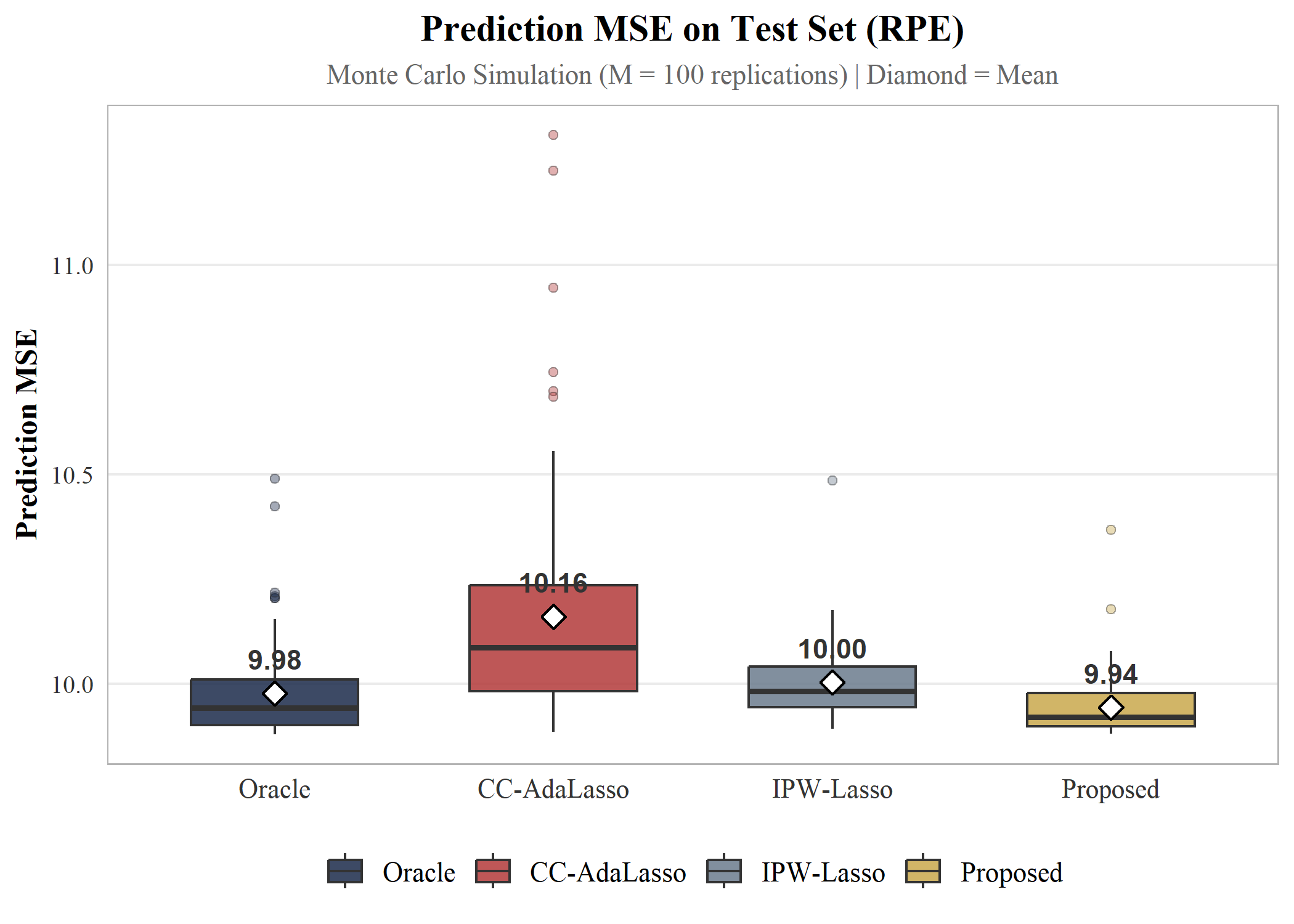}
\end{center}
\begin{center}
\underline{Figure 5.n=400  }
\end{center}

\begin{center}
\includegraphics[height=6cm]{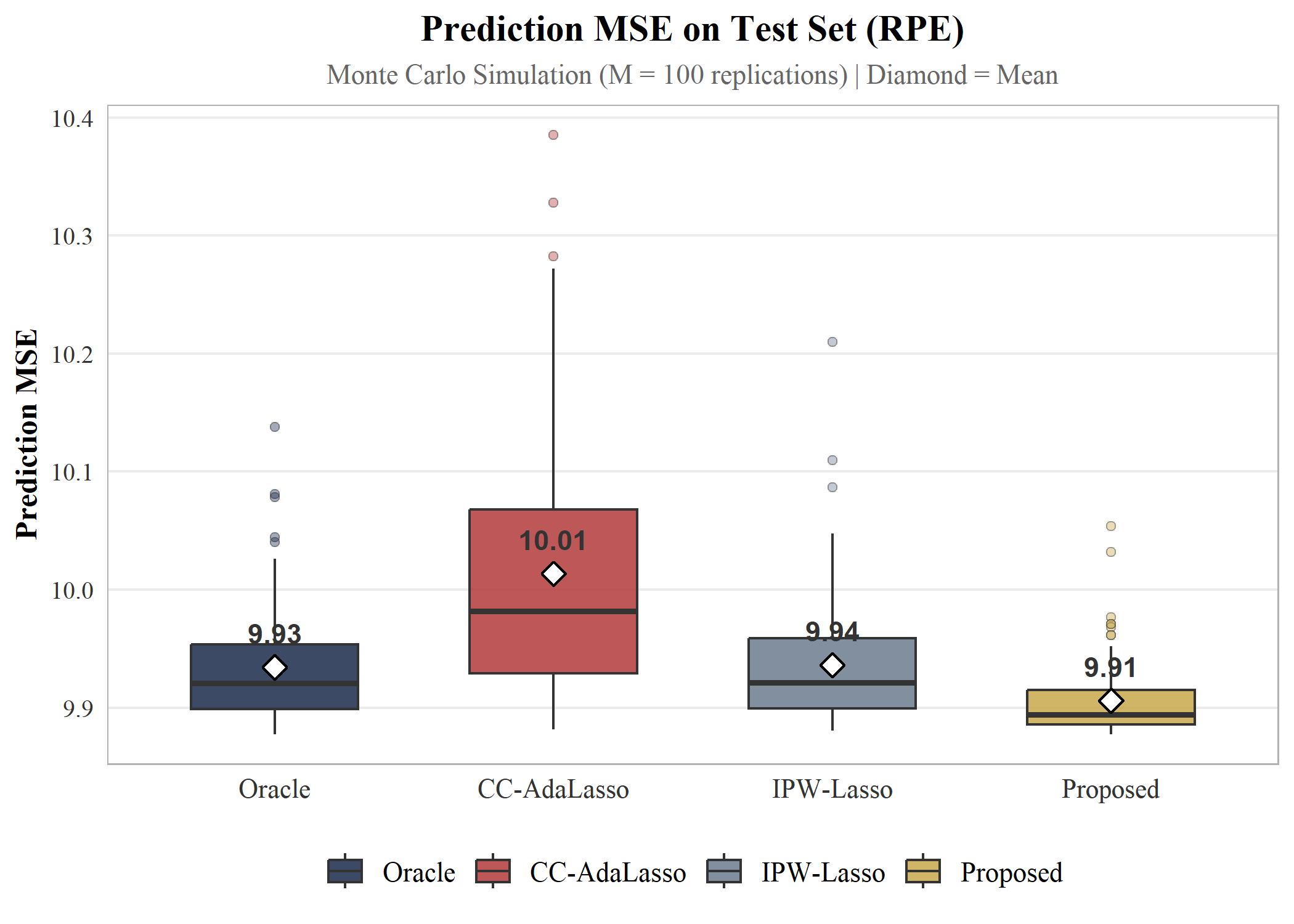}
\end{center}
\begin{center}
\underline{Figure 6.n=800  }
\end{center}

\begin{center}
\includegraphics[height=6cm]{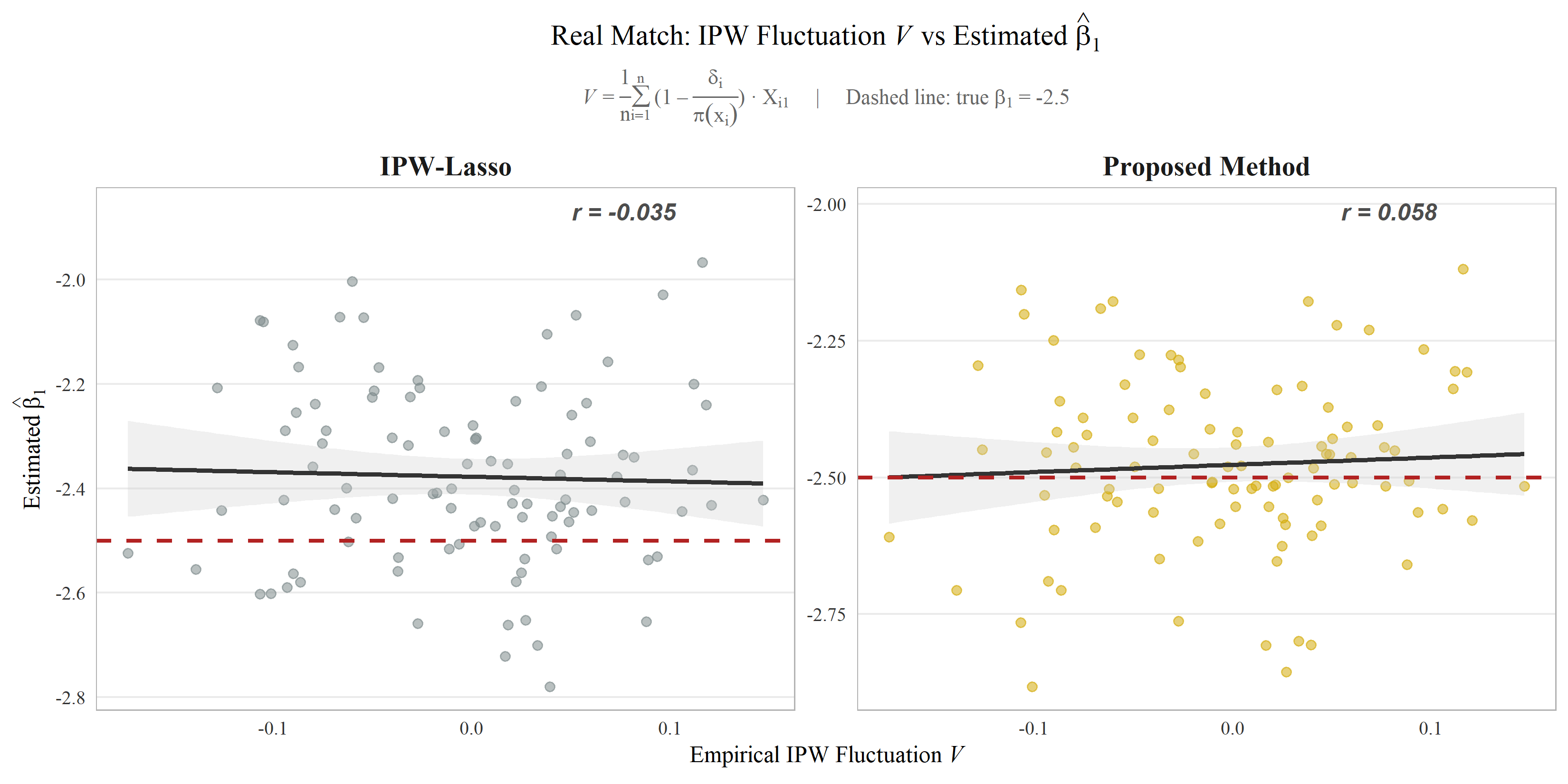}
\end{center}
\begin{center}
\underline{Figure 7.n=200  }
\end{center}

\begin{center}
\includegraphics[height=6cm]{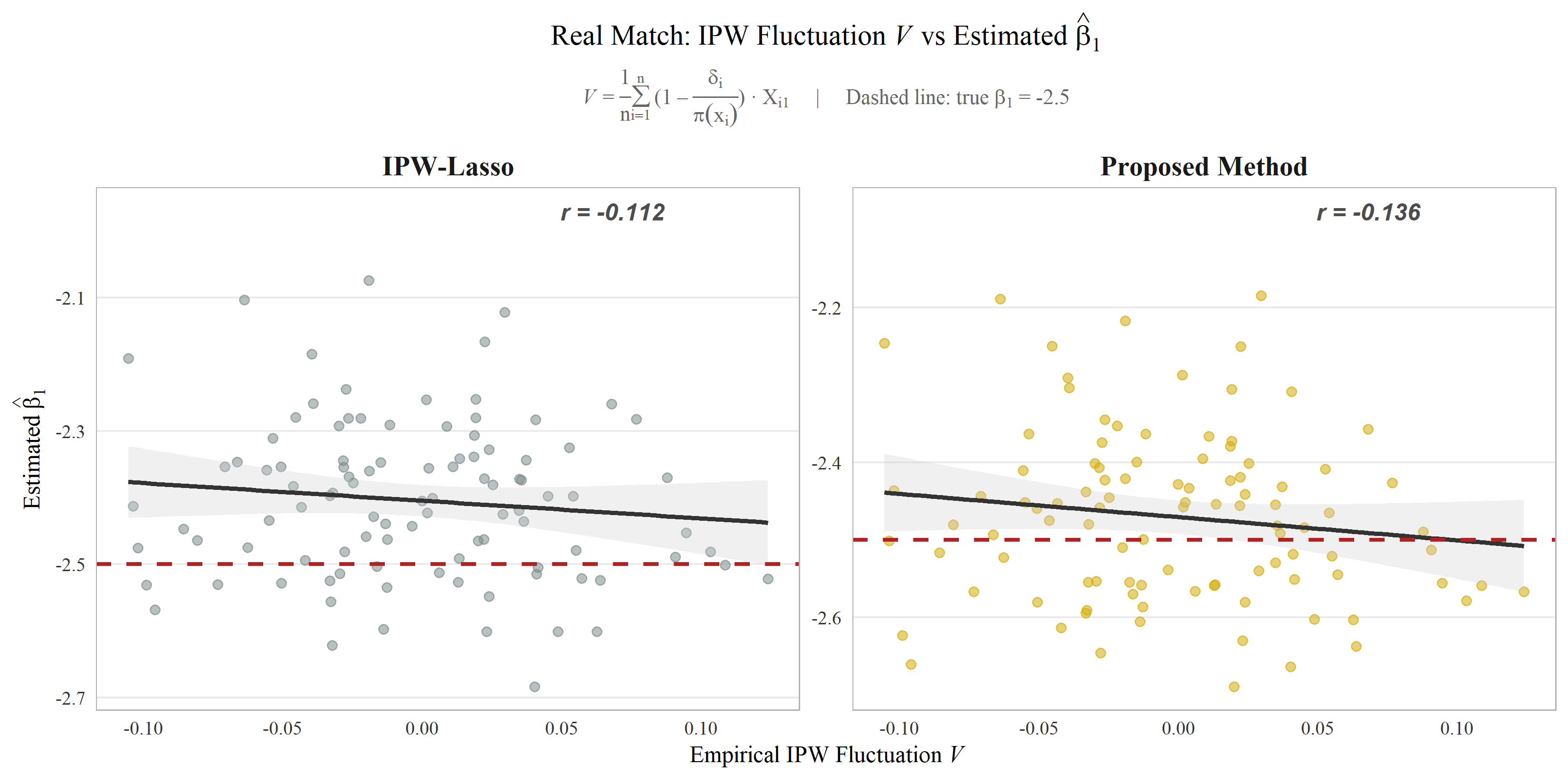}
\end{center}
\begin{center}
\underline{Figure 8.n=400  }
\end{center}

\begin{center}
\includegraphics[height=6cm]{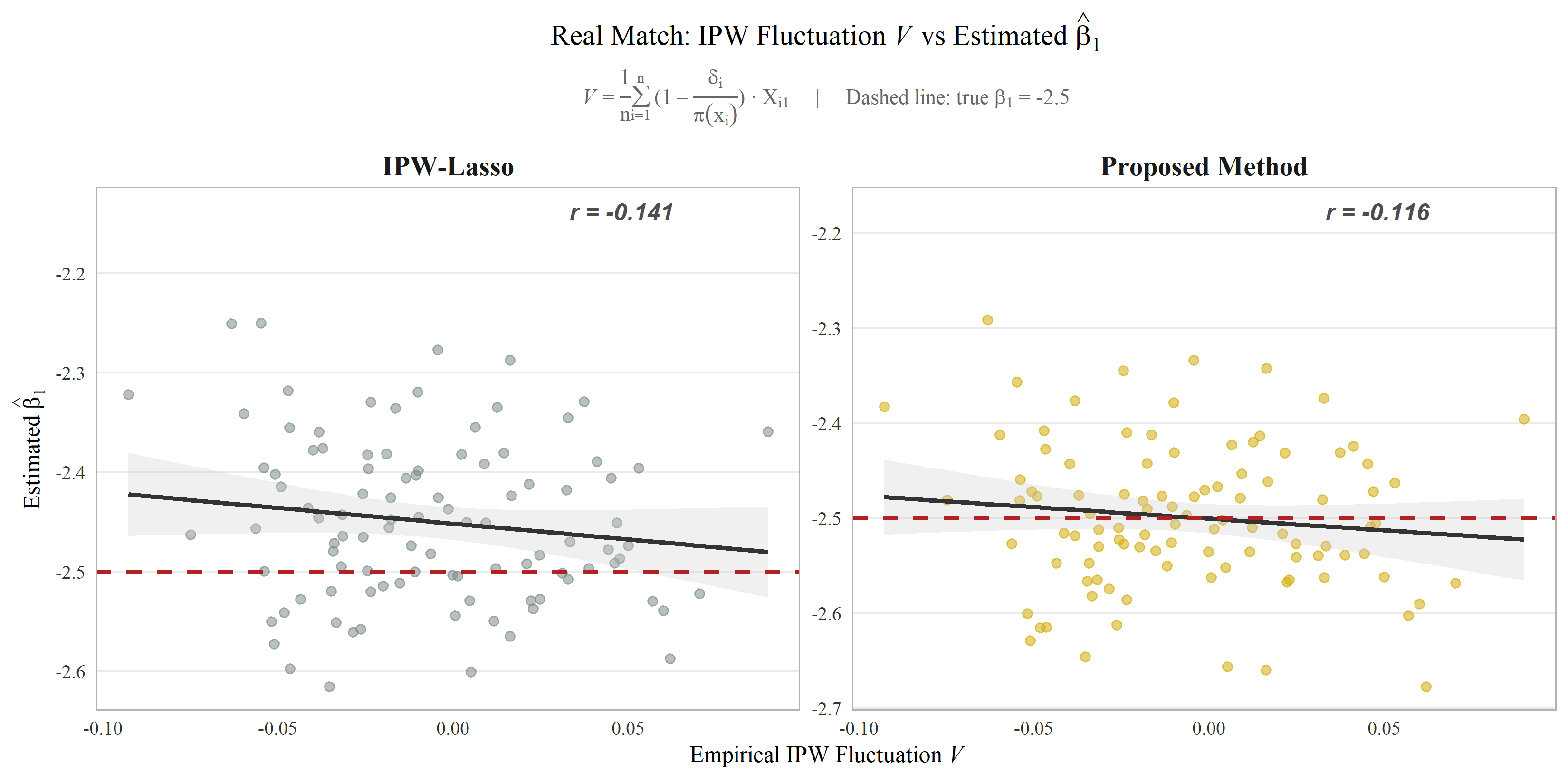}
\end{center}
\begin{center}
\underline{Figure 9.n=800  }
\end{center}

\begin{table}[ht]
\centering
\caption{Monte Carlo Simulation Results: Variable Selection and Prediction Performance} 
\label{tab:mc_results}
\scalebox{0.85}{
\begin{tabular}{lccccc}
  \toprule
Method & C Mean (SE) & I Mean (SE) & Exact Match & MSE ($\times 10^{-3}$) & RPE Mean (SE) \\ 
  \midrule
Oracle & 2.000 (0.000) & 0.000 (0.000) & 100.0\% & 9.67 (0.96) & 9.983 (0.010) \\ 
  CC-AdaLasso & 2.000 (0.000) & 0.890 (0.113) & 51.0\% & 32.50 (2.78) & 10.167 (0.023) \\ 
  IPW-Lasso & 2.000 (0.000) & 2.900 (0.190) & 11.0\% & 9.97 (0.70) & 10.011 (0.009) \\ 
  Proposed Method & 2.000 (0.000) & 0.800 (0.155) & 67.0\% & 5.84 (0.59) & 9.953 (0.007) \\ 
   \bottomrule
\end{tabular}
}
\end{table}

\begin{table}[ht]
\centering
\caption{Mean Coefficient Estimates vs True Values} 
\label{tab:coef_estimates}
\scalebox{0.85}{
\begin{tabular}{lcccccccc}
  \toprule
Method & Beta\_X1 & Beta\_X2 & Beta\_X3 & Beta\_X4 & Beta\_X5 & Beta\_X6 & Beta\_X7 & Beta\_X8 \\ 
  \midrule
True Value & -2.500 & 3.000 & 0.000 & 0.000 & 0.000 & 0.000 & 0.000 & 0.000 \\ 
  CC-AdaLasso & -2.434 & 2.919 & 0.009 & -0.013 & -0.013 & 0.019 & -0.002 & 0.002 \\ 
  IPW-Lasso & -2.404 & 2.889 & 0.003 & 0.000 & 0.005 & 0.014 & 0.003 & -0.002 \\ 
  Oracle & -2.490 & 2.974 & 0.000 & 0.000 & 0.000 & 0.000 & 0.000 & 0.000 \\ 
  Proposed Method & -2.470 & 2.952 & -0.001 & 0.004 & -0.001 & 0.008 & -0.000 & 0.003 \\ 
   \bottomrule
\end{tabular}
}
\end{table}

\section{Conclusion}
\qquad In this article, unlike existing predictive variables in the literature, we propose a relevant predictive variable framework by defining an adaptive Lasso with a weighted IPW factor and a KRR penalty function. We find that the weighted IPW lasso cannot be an oracle procedure while the adaptive lasso with weighted IPW factor and KRR penalty function can enjoy the oracle property. The numerical experimental results demonstrate that the MSE, RPE and RMR(the real match rate) obtained by our method are outperform competing methods in existing literature.. 

    In fact, experimental data demonstrate that the proposed adaptive weighted estimation method yields a substantially lower false selection rate of zero-coefficient terms compared with the conventional IPW-Lasso (reduced from 2.900 to 0.800). Its exact matching rate reaches $67.0/100$, considerably higher than the $11.0/100$ achieved by IPW-Lasso¹. Meanwhile, the proposed method attains a parameter estimation MSE of $5.84\times 10^{-3}$, which even outperforms the Benchmark Oracle estimator’s MSE of $9.67\times 10^{-3}$, demonstrating excellent near-minimax optimality.

\begin{center} {\bf \Large Appendix: Proofs} \end{center}

{\bf Proof of Lemma 1.} Define
$$
Z_n(\beta)=\frac{1}{n}[(Y-\sum_{j=1}^pX_j\beta_j)^TA(Y-\sum_{j=1}^pX_j\beta_j)]+\frac{\lambda_n}{n} \sum_{j=1}^p|\beta_j|^{\gamma}.
$$
If $\gamma\geq 1$, by the convexity of $Z_n(\beta)$, applying the standard results on the pointwise convergence in probability of $Z_n(\beta)$ to
$Z(\beta)+\sigma^2$(see \cite{Anderson1982} and \cite{Pollard1991}), we have
\begin{align}
\sup_{\beta\in K}|Z_n(\beta)-Z(\beta)-\sigma^2|\rightarrow_p 0\label{427s1}
\end{align}
for any compact set $K$ and
\begin{align}
\hat{\beta}^{(n)}=O_p(1).\label{427s2}
\end{align}
If $\gamma<1$, although $Z_n(\beta)$ is no longer convex, (\ref{427s1})
follows easily. Meanwhile, noticing
$$
Z_n(\beta)\geq\frac{1}{n}[(Y-\sum_{j=1}^pX_j\beta_j)^TA(Y-\sum_{j=1}^pX_j\beta_j)]
$$
and
$$
arg\min \frac{1}{n}[(Y-\sum_{j=1}^pX_j\beta_j)^TA(Y-\sum_{j=1}^pX_j\beta_j)]=O_p(1),
$$
we know that (\ref{427s2}) still holds.  \hfill$\Box$

{\bf Proof of Lemma 2.} Denoting $\epsilon=(\epsilon_1,...,\epsilon_n)^T$ and constructing
$$
V_n(u)=\left(\epsilon-\sum_{j=1}^pX_j\frac{u_j}{\sqrt{n}}\right)^TA\left(\epsilon-\sum_{j=1}^pX_j\frac{u_j}{\sqrt{n}}\right)-\epsilon^TA\epsilon
+\lambda_n\sum_{j=1}^p[|\beta_j+\frac{u_j}{\sqrt{n}}|^{\gamma}-|\beta_j|^{\gamma}],
$$
which is minimized at $\sqrt{n}(\hat{\beta}^{(n)}-\beta)$. Then it is easy to verify that
$$
\left(\epsilon-\sum_{j=1}^pX_j\frac{u_j}{\sqrt{n}}\right)^TA\left(\epsilon-\sum_{j=1}^pX_j\frac{u_j}{\sqrt{n}}\right)-\epsilon^TA\epsilon
\rightarrow_d -2u^Tw+u^TBu$$
with $w$ has a $N(0,\delta^2 B)$ distribution, and
$$
\lambda_n\sum_{j=1}^p[|\beta_j+\frac{u_j}{\sqrt{n}}|^{\gamma}-|\beta_j|^{\gamma}]\rightarrow\lambda_0\sum_{j=1}^pu_jsgn(\beta_j)|\beta_j|^{\gamma-1}\quad {\rm if }\ \gamma>1,
$$
$$
\lambda_n\sum_{j=1}^p[|\beta_j+\frac{u_j}{\sqrt{n}}|^{\gamma}-|\beta_j|^{\gamma}]\rightarrow\lambda_0\sum_{j=1}^p[u_jsgn(\beta_j)I(\beta_j\neq 0)+|u_j|I(\beta_j=0)]
$$
 if $\gamma=1$. By the finite-dimensional convergence, we know that $V_n(u)\rightarrow_d V(u)$.  Note that the convexity of $V_n$  and $V(u)$ has a unique minimum. We have
 $$
 arg\min V_n(u)=\sqrt{n}(\hat{\beta}^{(n)}-\beta)\rightarrow_d arg\min V(u).
 $$
Especially, $w$ has a $N(0, \sigma^2 B)$ distribution when $\lambda_0=0$. \hfill $\Box$

{\bf Proof Lemma 3.} Similar to Lemma 2, we construct
$$
V_n(u)=\left(\epsilon-\sum_{j=1}^pX_j\frac{u_j}{\sqrt{n}}\right)^TA\left(\epsilon-\sum_{j=1}^pX_j\frac{u_j}{\sqrt{n}}\right)-\epsilon^TA\epsilon
+\lambda_n\sum_{j=1}^p[|\beta_j+\frac{u_j}{\sqrt{n}}|^{\gamma}-|\beta_j|^{\gamma}].
$$
It is easy to see that
$$\lambda_n[|\beta_j+\frac{u_j}{\sqrt{n}}|^{\gamma}-|\beta_j|^{\gamma}]\rightarrow 0$$
if $\beta_j\neq 0$ because $\lambda_n=O(n^{\frac{\gamma}{2}})=o(\sqrt{n})$. Therefore,
$$
\lambda_n\sum_{j=1}^p[|\beta_j+\frac{u_j}{\sqrt{n}}|^{\gamma}-|\beta_j|^{\gamma}]\rightarrow\lambda_0\sum_{j=1}^p|u_j|^{\gamma}I(\beta_j=0)
$$
and consequently
$$
V_n(\cdot)\rightarrow_d V(\cdot)
$$
uniformly over $u$ in compact sets. Meanwhile, for all $u$ and $n$ large enough, there holds
\begin{align*}
V_n(u)&\geq \left(\epsilon-\sum_{j=1}^pX_j\frac{u_j}{\sqrt{n}}\right)^TA\left(\epsilon-\sum_{j=1}^pX_j\frac{u_j}{\sqrt{n}}\right)-\epsilon^TA\epsilon
-\lambda_n\sum_{j=1}^p|\frac{u_j}{\sqrt{n}}|^{\gamma}\\
&\geq \left(\epsilon-\sum_{j=1}^pX_j\frac{u_j}{\sqrt{n}}\right)^TA\left(\epsilon-\sum_{j=1}^pX_j\frac{u_j}{\sqrt{n}}\right)-\epsilon^TA\epsilon
-(\lambda_0+\delta)\sum_{j=1}^p|u_j|^{\gamma}\\
&:=V^{(l)}_n(u).
\end{align*}
Since $\gamma<1$, the quadratic terms of $V^{(l)}_n(u)$ grow faster
than the $|u_j|^{\gamma}$ terms, consequently, $arg\min V^{(l)}_n(u)=O_p(1)$ and $arg\min V_n(u)=O_p(1)$. Therefore,
$$arg\min V_n(u)\rightarrow_d arg\min V(u)$$
 because  $arg\min V(u)$ is unique with probability 1.\hfill$\Box$

{\bf Proof of Proposition 1.} Denote $u^*=arg\min V_2(u)$. Since $\sqrt{n}\hat{\beta}_{\mathcal{A}}\rightarrow_d u^*_{\mathcal{A}}$, we have
$\lim\sup_nP(\sqrt{n}\hat{\beta}_j=0,\forall j\notin \mathcal{A})\leq P(u^*_j=0,\forall j\notin \mathcal{A})$ due to the weak convergence result.
Meanwhile, $\mathcal{A}_n=\mathcal{A}$ implies that $\hat{\beta}_j=0$ for all $j\notin \mathcal{A}$, $P(\mathcal{A}_n=\mathcal{A})\leq P(\sqrt{n}\hat{\beta}_j=0, \forall j\notin \mathcal{A})$. Hence we only need to prove that $c=P(u^*_j=0, \forall j\notin \mathcal{A})<1$.
We discuss it in two cases.

Case 1. $\lambda_0=0$. In this case, it is easy to verify that $u^*=B^{-1}W\sim N(0,\sigma^2B^{-1})$ and $c=0$.

Case 2. $\lambda_0>0$. Since $v_2(u)$ is not differential at $u_j=0, \forall j \in \mathcal{A}$, we have
\begin{align}
&-2W_j+2(Bu^*)_j+\lambda_0 sgn(\beta^*_j)=0,\quad \forall j \in \mathcal{A},\label{57w1}\\
&|-2W_j+2(Bu^*)_j|\leq \lambda_0\label{57w2}
\end{align}
by the Karush-Kuhn-Tucker(KKT) optimality condition. Especially, if $u^*_j=0$ for all $j\notin \mathcal{A}$, then we can get
\begin{align}
-2W_{\mathcal{A}}+2B_{11}u^*_{\mathcal{A}}+\lambda_0 sgn(\beta^*_{\mathcal{A}})=0,\label{57w3}
\end{align}
and
\begin{align}
|-2W_{\mathcal{A}^c}j+2B_{11}u^*_{\mathcal{A}}|\leq \lambda_0\label{57w4}
\end{align}
componentwise by (\ref{57w1}) and (\ref{57w2}).  (\ref{57w3}) and (\ref{57w4}) imply that
$$
|-2W_{\mathcal{A}^c}+B_{21}B^{-1}_{11}(2W_{\mathcal{A}}-\lambda_0 sgn(\beta^*_{\mathcal{A}}))|\leq \lambda_0
$$
componentwise. Therefore,
$$
c\leq P(|-2W_{\mathcal{A}^c}+B_{21}B^{-1}_{11}(2W_{\mathcal{A}}-\lambda_0 sgn(\beta^*_{\mathcal{A}}))|\leq \lambda_0)<1. \qquad\qquad\qquad \qquad
 \Box
$$

{\bf Proof of Lemma 4:} Let $\beta=\beta^*+\frac{\lambda_n}{n}u$. Defining
$$
\Phi(u)=\left(\sum_{i=1}^n(1-\frac{\delta_i}{\pi(\tilde{x}_i)})(y_i-\sum_{j=1}^p(\beta^*_j+\frac{\lambda_n}{n}u_j)x_{ij})\right)^2
+\lambda_n\sum_{j=1}^p|\beta^*_j+\frac{\lambda_n}{n}u_j|,
$$
we can get
\begin{align*}
\Phi(u)-\Phi(0)&=-\left(\sum_{i=1}^n(1-\frac{\delta_i}{\pi(\tilde{x}_i)})2(y_i-\sum_{j=1}^p\beta^*_jx_{ij})\right)
\left(\sum_{i=1}^n(1-\frac{\delta_i}{\pi(\tilde{x}_i)})\sum_{j=1}^p\frac{\lambda_n}{n}u_jx_{ij}\right)\\
&\quad+\left(\sum_{i=1}^n(1-\frac{\delta_i}{\pi(\tilde{x}_i)})\sum_{j=1}^p\frac{\lambda_n}{n}u_jx_{ij}\right)^2
+\lambda_n\sum_{j=1}^p(|\beta^*_j+\frac{\lambda_n}{n}u_j|-|\beta^*_j|).
\end{align*}
Denote
\begin{align*}
V^{(n)}_3(u)&=u^T\left(\frac{X^TAX}{n}\right)u-2\frac{\epsilon^TAX}{\sqrt{n}}\frac{\sqrt{n}}{\lambda_n}u\\
&\quad+\sum_{j=1}^p\frac{n}{\lambda_n}(|\beta^*_j+\frac{\lambda_n}{n}u_j|-|\beta^*_j|).
\end{align*}
Then $\Phi(u)-\Phi(0)=\frac{\lambda^2_n}{n}V^{(n)}_3(u)$ and $$\hat{u}_n=arg\min \Phi(u)=arg\min V^{(n)}_3
=\frac{n}{\lambda_n}(\beta^{(n)}-\beta^*)$$ under the assumption $\frac{1}{n}X^TAX\rightarrow B$. Since $\frac{\epsilon^TAX}{\sqrt{n}}$ is
$O_p(1)$, and $\frac{\lambda_n}{\sqrt{n}}\rightarrow \infty$ equals to $\frac{\sqrt{n}}{\lambda_n}\rightarrow 0$, it is easy to obtain that
$\frac{\epsilon^TAX}{\sqrt{n}}\frac{\sqrt{n}}{\lambda_n}\rightarrow_p 0$ by Slutsky's theorem. Meanwhile,
\begin{align*}
\frac{n}{\lambda_n}(|\beta^*_j+\frac{\lambda_n}{n}u_j|-|\beta^*_j|)=\left\{\begin{array}{ll}
&u_jsgn(\beta_j)\quad {\rm if } \quad \beta_j\neq 0,\\
&|u_j|\qquad \qquad  {\rm if } \quad\beta_j=0.
\end{array}\right.
\end{align*}
Consequently, $V^{(n)}_3(u)\rightarrow_p V_n(u)$ for every $u$. Note that $V_3(u)$ has a unique minimizer due to $B$ is a positive definite matrix and
$V^{(n)}_3(u)$ is a convex function. We have $\hat{u}_n=arg\min (V^{(n)}_3)\rightarrow_p arg\min (V_3)$.\hfill $\Box$

{\bf Proof of Theorem 1:} We first prove the theorem when both the limit of $\frac{\lambda_n}{\sqrt{n}}$ and that of $\frac{\lambda_n}{n}$ exist. By the results
of Proposition 1, the lasso selection isn't consistent if $\frac{\lambda_n}{\sqrt{n}}\rightarrow \lambda_0\geq 0$. Therefore, under the assumption of the consistent selection, we only need to consider it in the three cases  below.

Case 1: $\frac{\lambda_n}{n}\rightarrow \infty$. In this case, $\hat{\beta}^{(n)}_j\rightarrow_p 0$ for all $j=1,2,...,p$ which is a contradiction to the assumption of the consistent selection.

Case 2: $\frac{\lambda_n}{n}\rightarrow \lambda_0$, $0<\lambda_0<\infty$. In this case, $\hat{\beta}^{(n)}\rightarrow_p \beta_*$ by Lemma 1, where $\beta_*$
is a nonrandom vector. And $\beta_{*j}=0$ for all $j\notin \mathcal{A}$ because $P(\mathcal{A}_n=\mathcal{A})=1$. Hence we will take a $j\in \mathcal{A}$ and
consider the event $j\in \mathcal{A}_n$.  First, we have
$$
-2x^T_jA(y-X{\hat{\beta}}^{(n)}+\lambda_n sgn ({\hat{\beta}}^{(n)}_j)=0
$$
by the KKT optimal conditions. Therefore,
$$P(j\in \mathcal{A}_n)\leq P(|-2\frac{x^T_jA(y-X{\hat{\beta}}^{(n)})}{n}|=\frac{\lambda_n}{n}.$$
Meanwhile,
$$
-2\frac{x^T_jA(y-X{\hat{\beta}}^{(n)})}{n}=-2\frac{x^T_jAX(\beta^*-{\hat{\beta}}^{(n)})}{n}-2\frac{x^T_jA\epsilon}{n}\rightarrow_p -2(B(\beta^*-\beta_*))_j.
$$
That is, $P(j\in \mathcal{A}_n)\rightarrow 1$ means that $|2(B(\beta^*-\beta_*))_j|=\lambda_0$. Concurrently, we can pick a $j'\notin \mathcal{A}$ which satisfies $P(j'\notin \mathcal{A}_n)\rightarrow 1$ and consider the event $j'\notin \mathcal{A}_n$. We get
$$
|-2x^T_{j'}A(y-X{\hat{\beta}}^{(n)}|\leq \lambda_n.
$$
Similarly,
$$P(j'\notin \mathcal{A}_n)\leq P(|-2\frac{x^T_{j'}A(y-X{\hat{\beta}}^{(n)})}{n}|\leq \frac{\lambda_n}{n}).$$
Consequently, $P(j'\notin \mathcal{A}_n)\rightarrow 1$ implies that $|2(B(\beta^*-\beta_*))_{j'}|\leq \lambda_0$.
Note that
$$
B(\beta^*-\beta_*)=\left(\begin{array}{cc}
         B_{11}(\beta^*_{\mathcal{A}}-\beta_{*\mathcal{A}})\\
          B_{21}(\beta^*_{\mathcal{A}}-\beta_{*\mathcal{A}}) \\
       \end{array}
\right)
$$
Denoted the sign vector of $B_{11}(\beta^*_{\mathcal{A}}-\beta_{*\mathcal{A}})$ by $s_*$. Then
\begin{align}
B_{11}(\beta^*_{\mathcal{A}}-\beta_{*\mathcal{A}})=\frac{\lambda_0}{2}s_*\label{941}
\end{align}
and
\begin{align}
|B_{21}(\beta^*_{\mathcal{A}}-\beta_{*\mathcal{A}})|\leq \frac{\lambda_0}{2}.\label{942}
\end{align}
 (\ref{941}) implies that $(\beta^*_{\mathcal{A}}-\beta_{*\mathcal{A}})=\frac{\lambda_0}{2}B_{11}^{-1}s_*$. Substituting it into (\ref{942}), we have
 \begin{align}
 |B_{21}B_{11}^{-1}s_*|\leq 1.\label{943}
 \end{align}

Case 3. In this case, $\frac{n}{\lambda_n}({\hat{\beta}}^{(n)})-\beta^*)\rightarrow_p u^*=arg \min(V_3)$ and $u^*$ is a nonrandom vector. Taking any $j\notin \mathcal{A}$, we have $u^*_j=0$ since $P({\hat{\beta}}^{(n)}_j)=0\rightarrow 1$ and $\frac{n}{\lambda_n}{\hat{\beta}}^{(n)}_j\rightarrow_p u^*_j$. Meanwhile, recall that
$$
V_3(u)=u^TBu+\sum_{j\in \mathcal{A}}[u_js^*_j]+\sum_{j\notin \mathcal{A}}|u_j|
$$
with $s^*=sgn (\beta^*_{\mathcal{A}})$. Similar to the discussions in Case 2, we have $u^*_{\mathcal{A}}=-\frac{1}{2}B^{-1}_{11}s^*$ and $u^*_{\mathcal{A}^c}=0$. Consider
$$
Z(r)=r^T(B_{22}-B_{21}B^{-1}_{11}B_{12})r-r^TB_{21}B^{-1}_{11}s^*+\sum_i|r_i|.
$$
It is easy to verify that $arg\min Z(r)=u^*_{\mathcal{A}^c}$. Then using KTT optimality condition, we can get
\begin{align}
 |B_{21}B_{11}^{-1}s_*|\leq 1.\label{951}
 \end{align}

For the general sequence of $\{\frac{\lambda_n}{\sqrt{n}}\}$ and $\frac{\lambda_n}{n}$, we can choose a subsequence $n_k$ such that
the limit of $\{\frac{\lambda_{n_k}}{\sqrt{n_k}}\}$ and $\{\frac{\lambda_{n_k}}{n_k}\}$ exist, consider the subsequence $\{{\hat{\beta}}^{*(n_k)}\}$
and establish the same conclusion.\hfill $\Box$

$$
u^*_{\mathcal{A}}=-B^{-1}_{11}(B_{12}u^*_{\mathcal{A}^c}+\frac{1}{2}s^*).
$$

$$
B=\left(
       \begin{array}{cc}
         B_{11} & B_{12}\\
          B_{21} & B_{22} \\
       \end{array}
     \right)
$$

{\bf The proof of Corollary 1:} Computing directly, we can obtain
$$
B^{-1}_{11}=\frac{1}{1-\rho_1}[I_1-\frac{\rho_1}{[(n+1)-\rho_1]J_1}]
$$
and
\begin{align*}
B_{21}B^{-1}_{11}=\frac{(n+1)\rho_2\vec{1}}{(n+1)-\rho_1}.
\end{align*}
Hence
$$
B_{21}B^{-1}_{11}s=\frac{(n+1)\rho_2}{(n+1)-\rho_1}\left(\sum_j^{p_0}s_j \right).
$$
 Note that $|\sum_j^{p_0}s_j|\geq 1$ when $p_0$ is odd. Then $|B_{21}B^{-1}_{11}s|\geq 1$ if $\frac{(n+1)\rho_2}{(n+1)-\rho_1}>1$.
 \hfill $\Box$

{\bf The proof of Theorem 2:} Letting $\beta=\beta^*+\frac{u}{\sqrt{n}}$ and considering
$$
\Psi_n(u)=[Y-\sum_{j=1}^px_j(\beta^*_j+\frac{u_j}{\sqrt{n}})]^TA[Y-\sum_{j=1}^px_j(\beta^*_j+\frac{u_j}{\sqrt{n}})]
+\lambda_n\sum_{j=1}^p\hat{w}_j|\beta^*_j+\frac{u_j}{\sqrt{n}}|,
$$
we first hope to prove the part of asymptotic normality. Denote $\hat{u}^{(n)}=arg \min \Psi_n(u)$. We have
$\hat{\beta}^{*(n)}=\beta^*+\frac{\hat{u}^{(n)}}{\sqrt{n}}$, or $\hat{u}^{(n)}=\sqrt{n}(\hat{\beta}^{*(n)}-\beta^*)$.
Let
\begin{align*}
V^{(n)}_4(u)\equiv u^T\left(\frac{1}{n}X^TAX\right)u-2\frac{\epsilon^TAX}{\sqrt{n}}u
+\frac{\lambda_n}{\sqrt{n}}\sum_{j=1}^p\hat{w}_j\sqrt{n}(|\beta^*_j+\frac{u_j}{\sqrt{n}}|-|\beta^*_j|).
\end{align*}
Obviously, $V^{(n)}_4(u)=\Psi_n(u)-\Psi_n(\bf{0})$. Recall that $\frac{1}{n}X^TAX\rightarrow B$ and $\frac{\epsilon^TAX}{\sqrt{n}}u\rightarrow_d W=N(\bf{0},\sigma^2 B)$. Now we consider the term $\frac{\lambda_n}{\sqrt{n}}\sum_{j=1}^p\hat{w}_j\sqrt{n}(|\beta^*_j+\frac{u_j}{\sqrt{n}}|-|\beta^*_j|)$.

If $\beta^*=0$, then $\sqrt{n}(|\beta^*_j+\frac{u_j}{\sqrt{n}})=|u_j|$ and $\frac{\lambda_n}{\sqrt{n}}\hat{w}_j=\frac{\lambda_n}{\sqrt{n}}n^{\frac{\gamma}{2}}
(|\sqrt{n}\hat{\beta}_j|)^{-\gamma}$ with $\sqrt{n}\hat{\beta}_j=O_p(1)$.

If $\beta^*\neq 0$, note that $\hat{w}_j\rightarrow_p |\beta^*_j|^{-\gamma}$ and $\sqrt{n}(|\beta^*_j+\frac{u_j}{\sqrt{n}}|-|\beta^*_j|)\rightarrow u_j sgn(\beta^*_j)$, we have $\frac{\lambda_n}{\sqrt{n}}\hat{w}_j\sqrt{n}(|\beta^*_j+\frac{u_j}{\sqrt{n}}|-|\beta^*_j|)\rightarrow 0$ by Slutsky theorem.

Applying Slutsky theorem again, we obtain $V^{(n)}_4(u)\rightarrow_d V_4(u)$ for every $u$, where
$$
V_4(u)=\left\{\begin{array}{cc}
        u^T_{\mathcal{A}}B_{11}u_{\mathcal{A}}-2u^T_{\mathcal{A}}W_{\mathcal{A}}& {\rm if}\  u_j=0,\ \forall j\notin \mathcal{A},\cr\\
          \infty & {\rm otherwise} \\
       \end{array}\right.
$$
$V^{(n)}_4$ is convex and attaches its minimum at $(B^{-1}_{11}W_{\mathcal{A}},\bf{0})^T$. Recall that $W=N(\bf{0},\sigma^2 B)$, hence $W_{\mathcal{A}}=N(\bf{0}, \sigma^2 B_{11})$. Finally, by the epi-convergence results of Geyer(1994) as well as those of Knight and Fu(2000), we get
\begin{align}
u^{(n)}_{\mathcal{A}}\rightarrow_d B^{-1}_{11}W_{\mathcal{A}}\quad {\rm and}\quad u^{(n)}_{\mathcal{A}^c}\rightarrow_d 0,\label{952}
\end{align}
which finishes the proof of asymptotic normality.

We will prove the consistency below. Using the property of asymptotic normality, we have ${\hat{\beta}}^{(n)}_j\rightarrow_p \beta^*_j$ and $P(j\in \mathcal{A}^*_n)\rightarrow 1$. So we only need to show that $P(j'\notin \mathcal{A}^*_n)\rightarrow 0$ for any $j'\notin \mathcal{A}^*_n$. Hence we consider the event $j'\in \mathcal{A}^*_n$. The KTT conditions implies that $2x^T_{j'}A(y-X\beta^{*(n)})=\lambda_n\hat{w}_{j'}$.
Noticing the following two facts
\begin{align*}
&\frac{\lambda_n\hat{w}_j}{\sqrt{n}}=\frac{\lambda_n}{\sqrt{n}}n^{\frac{\gamma}{2}}\frac{1}{|\sqrt{n}\hat{\beta}_{j'}|^{\gamma}}\rightarrow_p \infty,\\
&2\frac{x^T_{j'}A(y-X{\hat{\beta}}^{*(n)})}{\sqrt{n}}=2\frac{x^T_{j'}AX\sqrt{n}(\beta^*-{\hat{\beta}})}{n}+2\frac{x^T_{j'}A\epsilon}{\sqrt{n}},
\end{align*}
it can be verified that
$$
2\frac{x^T_{j'}AX\sqrt{n}(\beta^*-{\hat{\beta}})}{n}\rightarrow_d \ {\rm some \ normal  \ distribution}
$$
and
$$
2\frac{x^T_{j'}A\epsilon}{\sqrt{n}}\rightarrow_d\quad  N(\bf{0},4\|x_{j'}\|^2\sigma^2).
$$
Therefore,
$$
P(j'\in \mathcal{A}_n^{\infty})\leq P(2x^T_{j'}A(y-X{\hat{\beta}}^{*(n)})=\lambda_n\hat{w}_{j'})\rightarrow 0.
$$
Theorem 2 is proved. \hfill $\Box$

\end{document}